\documentstyle[icmp]{article}

\def\Journal#1#2#3#4{{#1} {\bf #2}, #3 (#4)}

\def\JMP{\em J. Math. Phys.}

\def\JPA{\em J. Phys. A, Math. Gen.}

\def\be{\begin{equation}}
\def\ee{\end{equation}}
\def\bea{\begin{eqnarray}}
\def\eea{\end{eqnarray}}
\def\nn{\nonumber \\}
\def\openone{\mbox{1\kern -0.25em I}}
\def\openR{\mbox{I\kern-0.20em R}}
\def\con{\,\rule{1ex}{0.75pt}\rule{0.75pt}{1ex}\,}
\def\Cl{C\kern -0.2em\ell}
\def\us{{\mbox{\underline{$s$}}}}
\def\ub{{\mbox{\underline{$b$}}}}
\def\ob{{\mbox{$\overline{b}^\epsilon$}}}

\begin{document}

\title{ON THE RELATION OF CLIFFORD-LIPSCHITZ GROUPS
       TO $q$-SYMMETRIC GROUPS}

\author{BERTFRIED FAUSER}

\address{Universit\"at Konstanz,
Fakult\"at f\"ur Physik,
Fach M 678\\
78457 Konstanz, Germany\\
E-mail: Bertfried.Fauser@uni-konstanz.de} 


\maketitle\abstracts{
It can be shown that it is possible to find a representation of
Hecke algebras within Clifford algebras of multivectors. These
Clifford algebras possess a {\em unique}\/ gradation and a
possibly {\em non-symmetric}\/ bilinear form. Hecke algebra
representations can be classified, for non-generic $q$, by
Young tableaux of the symmetric group due to the isomorphy of the 
group algebras for $q \rightarrow 1$. Since spinors can be constructed 
as elements of minimal left (right) ideals obtained by the left (right) 
action on primitive idempotents, we are able to construct $q$-spinors
from $q$-Young operators corresponding to the appropriate symmetry
type. It turns out that an anti-symmetric part in the Clifford bilinear 
form is necessary. $q$-deformed reflections (Hecke 
generators) can be obtained only for {\em even}\/ multivector 
aggregates rendering this symmetry a composite one. In this 
construction one is able to deform spin groups only, though not 
pin groups. The method is closely related to a projective 
interpretation.}

\section{Introduction}

We recently showed that it is possible to find linear representations
of Hecke algebras within Clifford algebras of multivectors.\cite{Hecke1}
These Clifford algebras are characterized by an arbitrary bilinear form,
while a quadratic form does possess, in the case of characteristic zero,
a symmetric bilinear form only. The aim of this note is to show in which
way $q$-spinor representations and Hecke algebra representations are 
united due to $q$-deformed Young operators. $q$-deformed spinors and 
$q$-spin groups result.

Starting from non-deformed algebras and embedding the
deformed ones, we contrast recent work on $q$-spinors based on 
$q$-deformed Clifford algebras.\cite{Bautista,Oziewicz} 
Our method is conservative in the sense that the usual interpretation
of the undeformed Clifford algebra remains valid, while the $q$-deformed
objects exhibit a novel composite structure. Computations were performed 
by {\em Clifford},\/ a Maple V add-on.\cite{Ablamowicz}

\section{Hecke algebra representations within Clifford algebras
of multivectors}

Let $\bigwedge (V)$ be the Grassmann algebra built over the linear space
$V$ of dimension $2n$. We denote the linear space underlying $\bigwedge
(V)$ as $F\!\bigwedge (V)$. A Clifford map $\gamma$ is a linear map from
$V$ into a unital associative algebra ${\cal A}$, which satisfies 
\bea\label{CliffordMap}
&\gamma_x^2 \,=\, \gamma_x \gamma_x \,=\, Q(x)\openone.&
\eea
This can be polarized to $\gamma_x \gamma_y + \gamma_y \gamma_x = 2G(x,y)=Q(x+y)-Q(x)-Q(y)$, 
where $G$ is the symmetric bilinear form associated to the quadratic
form $Q$. We observe with Chevalley\cite{Chevalley} that for $x \in V$
a linear combination $\gamma_x := x\con + x\wedge$ is a Clifford map
operating on $V$. For details see\cite{thesis,Chap4,Hecke1}. Lifting
this map to $F \bigwedge (V)$ using
\bea
i)  && x \con y \,=\, < x \mid y > \,=\, B(x,y) \nn
ii) && x \con (u \wedge v) \,=\, (x \con u) \wedge v
                              + (\hat{u}) \wedge (x \con v) \nn
iii)&& (u \wedge v) \con w \,=\, u \con ( v \con w),
\eea
where $x,y \in V$, $u,v,w \in \bigwedge (V)$ and $\hat{\ }$ is the
main involution $\hat{V}=-V$, results in the Clifford algebra of
multivectors. Note, $B$ is an arbitrary bilinear form, however
the anticommutator relation is not altered by an
antisymmetric part of the bilinear form $B$. The corresponding
algebras are isomorphic as Clifford algebras, not however, as graded
algebras.

By choosing a basis $\{e_i\}$ and a suitable bilinear form:
(appropriate indices only, $i,j \in \{1,\cdots,2n\}$)
\bea
[B_{ij}]=
[B(e_i,e_j)]&:=&q[\delta_{i,j-n}]+
\lambda[\delta_{i-n,j-1}]+[\delta_{i-n,j}]
+\frac{q}{\lambda}[\delta_{i-n-1,j}],
\eea
we are able to find $n$ bivector elements $b_i:= e_i \wedge e_{n+i}$
satisfying the Hecke relations
\bea\label{Hecke-rel}
i)  && b_i^2 \,=\, (1-q)b_i +q \nn
ii) && b_i b_j \,=\, b_j b_i \quad \vert i-j\vert \ge 2 \nn
iii)&& b_i b_{i+1} b_i \,=\, b_{i+1} b_i b_{i+1} 
\eea
thereby providing a Hecke algebra representation.\cite{Hecke1}
 
\section{$q$-spinor representations from $q$-Young operators}

Spinors can be defined as elements of minimal left (right)
ideals of Clifford algebras. Iff $f$ is a primitive (irreducible) 
idempotent, one has ${\cal S}= \Cl f$, $f^2=f$. A primitive idempotent
element might be written in the form $f=\prod 1/2(1\pm x_i)$ with
$[x_i,x_j]_-=0$, $x_i^2=1$ and a choice of signs. From $i)$ in
(\ref{Hecke-rel}) we obtain mutually annihilating projectors due to $(q+b_i)(1-b_i)=0$,
hence
\bea
P^{q+}_i:=\frac{1}{1+q}(q+b_i)
&\quad& 
P^{q-}_i:=\frac{1}{1+q}(1-b_i).
\eea
However, $P^{q\pm}_i P^{q\pm}_j$ do not necessarily commute iff 
$i=j\pm1$. 

We will specialize to $H_3(q)$ resp. $S_3$ in the sequel, so
$n=2$ and $\Cl (V,B)\cong \Cl_{2,2}(\openR,B)$. An algebraic basis 
of $H_3(q)$ is given by the elements 
$\{\openone, b_1, b_2, b_{12},$ $ b_{21}, b_{121} \}$, 
using the abbreviations $b_{ij}:=b_ib_j$ etc.
These are only some elements of the even subalgebra 
$\Cl^+_{2,2}$ of $\Cl(\openR^4,B)$. But $e_{1\wedge 2}=
e_1\wedge e_2$ and $e_{3\wedge 4}$ can not be built from the $b_i$s.
We thus have $H_3(q) \subset \Cl^+(\openR^4,B)$. This is however
not true for higher dimensions.\cite{Hecke1}

Composing the row symmetrizer $R^q$ and column antisymmetrizer $C^q$
we can construct the following Young operators $Y_k=C^q_{k^\prime}
R^q_{k^{\prime\prime}}$ in accordance with\cite{q-young-op}
\bea
Y_{sym}&:=&\frac{1}{q^2+q+1}\Big\{
q^2\openone+q(e_{1\wedge 3}+e_{2\wedge 4})
-e_{1\wedge 2\wedge 3\wedge 4}\Big\} \nn
Y{{12} \atop 3} &:=& \frac{1}{(q+1)(q^2+q+1)}\Big\{
q\openone+e_{1\wedge 3}-q(q+1)e_{2\wedge 4}+
(q+1)e_{1\wedge 2\wedge 3\wedge 4}\Big\} \nn
Y{{13} \atop 2} &:=& \frac{1}{(q+1)(q^2+q+1)}\Big\{
q(2q+1)\openone-q^2 e_{1\wedge 3}
+(q+1) e_{2\wedge 4}
\nn &&
-\lambda (q+1) e_{1\wedge 4}
-q^2 e_{1\wedge 3}
+(q+1) e_{1\wedge 2\wedge 3\wedge 4}\Big\} \\
Y_{asym} &:=&\frac{1}{q^2+q+1}\Big\{
(1-q)\openone-e_{1\wedge 3}-e_{2\wedge 4}
+\lambda e_{1 \wedge 4}+\frac{q}{\lambda}e_{2\wedge 3}
-e_{1\wedge 2\wedge 3\wedge 4}\Big\} .\nonumber
\eea
They are normalized, mutually annihilating $Y_l Y_k = \delta_{ik}Y_k$ 
and complete: $\openone = Y_{sym}+Y{{12}\atop 3}
+Y{{13}\atop 2}+Y_{asym}$. 
In the $H_3 (q)$ case $q$ has to be unequal to $exp(2\pi i/n)$ 
for $n\in\{1,2,3\}$: non generic $q$.
These Young operators provide $q$-idempotents for the construction 
of left ideals, which carry a representation of $H_3(q)$. 
The dimension of such left ideals follows from the degeneracy of the
eigenvalues of the transposition class-sum $C_n$ e.g. 
$C_3:=b_1+b_2+1/q\, b_{121}$ 
which is a central operator and characterizes uniquely all 
representations of $H_n(q)$.\cite{transposition class-sum} 
We obtain two one-dimensional representations, the roots of $i)$ 
in (\ref{Hecke-rel}) belonging to $Y_{sym}$ and $Y_{asym}$, 
while $Y{{12}\atop 3}$ and $Y{{13}\atop 2}$ each generate
a two-dimensional space. The leftregular representation is
the direct sum of these spaces
\bea
{\cal S}_{reg} &:=& Y_{sym}
\oplus \left( {Y{{12}\atop 3} \atop b_2Y{{12}\atop 3}} \right)
\oplus \left( {Y{{13}\atop 2} \atop b_2Y{{13}\atop 2}} \right) 
\oplus Y_{asym}.
\eea
However, for the whole, even even part of the Clifford algebra 
this is not a faithful spinor representation, since 
$e_{1\wedge 2}\,{\cal S}_{reg} = 0 = e_{3\wedge 4}\,{\cal S}_{reg}$ 
annihilate this space.

Dealing with the even elements only, we have to note
that $\Cl^+_{2,2} \cong \Cl_{1,2}$ has a real 4-dimensional
spinor space: $\Cl_{1,2} \cong \openR(4)$. However, our mixed
symmetry operators can not be added to such a space, but we can
construct it by adjoining an odd element: $u=e_1+e_3$, 
$u^2=1+q$.
\bea
{\cal S}{{12}\atop 3} &:=& 
\left( {Y{{12}\atop 3} \atop b_2Y{{12}\atop 3}} \right)
\oplus
\left( {Y{{12}\atop 3} \atop b_2Y{{12}\atop 3}} \right) u
\eea
represents all 6 bi-vector elements and thus $\Cl^+_{2,2} 
\cong \Cl_{1,2}$. We have thus succeeded in finding a 
$q$-deformed mixed ${{12}\atop 3}$-symmetry type spinor in the 
even sector (spinor sector) of $\Cl_{2,2}$. The ${{13}\atop 2}$-type
is analogous.

\section{$q$-reflections and $q$-spin groups}

A theorem of Cartan states, that $n$-dimensional pseudo 
orthogonal groups can be constructed by less or equal than $n$ 
reflections $\us_i$. The action of such a reflection might be 
defined as a conjugation $\us_i(x) = s_i x s_i^{-1}$ or on spinors
as left translation $\us_i(\psi):= s_i \psi$. Furthermore, reflections 
satisfy the relations (\ref{Hecke-rel}) for $q\rightarrow 1$. 
Taking this into consideration one has $\us_i$ involutory: 
$\us_i \us_i=Id$ and $\us_i^{-1}= \us_i$. Now the $\ub_i$ might
serve as $q$-reflections in a $q$-analogue of CartanÂs theorem. 
However, since the $\ub_i$ are not involutions some care is needed.

Fortunately a neat definition of an (anit)involution on the
algebra is given by Clifford algebraic considerations. With the new 
main involution $(V\wedge V)\hat{~}=-V\wedge V$ 
and the reversion (of products) $(AB)\tilde{~}=\tilde{A}\tilde{B}$,
we compose the antiinvolution $\alpha_\epsilon$ as 
$\alpha_\epsilon(AB) =\alpha_\epsilon(B)\alpha_\epsilon(A)$, 
where $\epsilon=\pm 1$ incorporates the main involution or not. 
Using the reversion of the $e_i$--Clifford elements, we obtain:
\bea
&\ob_i \, :=\, \epsilon \tilde{b}_i \,=\,
\epsilon[(1-q)-b_i] \nn
&b_i\ob_i \,=\, \epsilon b_i[(1-q)-b_i]\,=\, -\epsilon q
\quad
&\Rightarrow
\quad
b_i^{-1} \,=\, -\frac{\epsilon}{q} \ob_i.
\eea
Such conjugations are connected to special elements and the
structure of Clifford algebras in a natural way.\cite{Hahn}
Exactly this conjugation interchanges the roots of i) in 
(\ref{Hecke-rel}) suggesting the choice $\epsilon=-1$. Furthermore,
this Clifford algebraic antiinvolution interchanges the row and 
column-parts in the Young operators $\alpha_\epsilon(C^q_k)\cong R^q_k$, 
$\alpha_\epsilon(R^q_k)\cong C^q_k$, as is easily seen in 
$\alpha_\epsilon(P_i^{q\pm})=\epsilon P_i^{q\mp}$. 
The Hecke-versor subgroup is thus given by 
$\underline{a}(x)= a x a^{-1}$, or $\underline{a}(\psi) = a \psi$
on spinors iff $a^{-1}:= a/(a\overline{a}^\epsilon)$ exists. 
The whole Clifford-Lipschitz or versor group
might result by adding an odd element once again, which now 
however has to be central and will be the {\em special}\/ element 
mentioned above.

To see that products of the $\ub_i$ constitute an invariance 
group of an inner product, we define the adjoint of 
$A=b_{i_1}\cdots b_{i_n}$ to be
$A^* \,:=\, \overline{A}^\epsilon = \epsilon^n 
(\tilde{b}_{i_n} \cdots \tilde{b}_{i_1})$, and for the spinor space
as
${\cal S}^*_q \,:=\, \overline{{\cal S}}^\epsilon_q$.
The Clifford--Lipschitz group of versors can then be defined with
$\Phi^q_\epsilon(x,y):=\overline{x}^\epsilon y$ as
\bea
\Gamma^{q+}& := &
\{ x \mid x\in \Cl^+(V,B),\quad \Phi^q_\epsilon(x,x)=1 \},
\eea
which reduces to the spin group for $q \rightarrow 1$. Furthermore, 
one can relate to such invariance groups an inner product, as the 
restriction of $\Phi^q_\epsilon$ to ${\cal S}_q$
$\Phi^q_\epsilon(x,y)\vert_{{\cal S}_q} \cong \phi^q_\epsilon
\cong {\cal S}^*_q {\cal S}_q$ shows, see\cite{Lounesto} for $q=1$.
This completes the claim that we have indeed constructed $q$-spin
groups and that our ideals contain $q$-spinors. The 
$\ub_i$-reflections generate the $q$-Weyl group. From 
(\ref{CliffordMap}) we see that $q$-pin groups, odd 
versor groups, cannot be constructed in this way.

The projective interpretation follows from the
{\em projective split}\/ defined by Hestenes,\cite{HZ,H} 
which allows the bi-vectors to be the projective points 
(line geometry). This is also reflected in the $2 \times 2$
block-matrix structure of $B$.

\section*{Acknowledgments}
My gratitude is expressed to the organizers of the ICGTMP98 conference
for financial support and to Kirsten Magee for improving my English.

\end{document}